\documentstyle[12pt]{article}
\title{ CLASSICAL YANG- BAXTER EQUATION AND LOW DIMENSIONAL TRIANGULAR
LIE BIALGEBRAS OVER ARBITRARY FIELD
 \thanks {This work is supported by National Science Foundation (No:19971074)}}
\author{Shouchuan Zhang \\ Department  of Mathematics, Hunan Normal   University,
410081 \\ P.R.China}
\date{}
\begin{document}
\newtheorem{Theorem}{\quad Theorem}[section]
\newtheorem{Proposition}[Theorem]{\quad Proposition}
\newtheorem{Definition}[Theorem]{\quad Definition}
\newtheorem{Corollary}[Theorem]{\quad Corollary}
\newtheorem{Lemma}[Theorem]{\quad Lemma}
\newtheorem{Example}[Theorem]{\quad Example}
\maketitle
\begin {abstract}
Let  $ L $  be a  Lie  algebra over  arbitrary  field   $ k $   with
dim  $ L $ =3  and  dim   $ L' $ =2. All   solutions  of  constant
classical Yang- Baxter  equation (CYBE) in  Lie  algebra   $ L $ are
obtained and  the  necessary  conditions  which   $ (L,[\
],\Delta_{r}, r) $  is a   coboundary  (  or   triangular  )  Lie
bialgebra    are given.


Key  words : Lie  bialgebra; Yang-Baxter  equation

\end {abstract}
\noindent \addtocounter{section}{0}

\section {Introduction}

The  concept  and  structures  of  Lie  coalgebras  were introduced
and  studied  by  W. Michaelis  in  \cite {Mi80,  Mi85}.
V.G.Drinfel'd and A. A. Belavin  in \cite {BD82,  Dr86}   introduced
the notion of triangular,  coboundary   $ L $   associated  to  a
solution   $ r\in L\otimes  L $   of the CYBE  and  gave  a
classification  of solutions of  CYBE with parameter for simple  Lie
algebras. W.Michaelis  in \cite {Mi94} obtained the structure  of  a
triangular,  coboundary Lie bialgebra  on any Lie  algebra
containing  linearly independent elements   $ a $  and  $ b $
satisfying  $ [a, b] = \alpha b $   for some non-zero  $ \alpha  \in
k $ by setting  $ r=a \otimes  b - b \otimes a $.

The Yang-Baxter equation  first came up in a paper by Yang as
factorition   condition of  the  scattering  S-matrix in the
many-body problem in one dimension and in work of Baxter on exactly
solvable models in statistical mechanics. It has been playing an
important role in mathematics and physics ( see \cite {BD82,
YG89}).Attempts to find solutions The Yang- Baxter equation in a
systematic way have let to the theory of quantum groups. The Yang-
Baxter  equation is of many forms. The classical Yang- Baxter
equation is one.

In many applications one need to know the solutions of classical
Yang- Baxter equation and know if a Lie algebra is a coboundary Lie
bialgebra or a triangular Lie bialgebras. A systematic study of low
dimensional Lie algebras, specially, of those Lie algebras that play
a role in physics( as e.g. sl(2, C),  or the Heisenberg algebra ),
is very useful.The author \cite {Zh98b,  Zh99f} obtained all
solutions of constant classical Yang- Baxter equation (CYBE) in Lie
algebra  $ L $  and give the sufficient and necessary conditions
which ( $ L, [ \ ], \Delta_{r},  r $ )is a coboundary ( or
triangular) Lie bialgebra with dim  $ L\leq 3 $  except the below
case : $ L $  is a Lie algebra over arbitrary field  $ k $  with dim
 $ L $ =3 and dim  $ L^{'}=2 $  .We shall resolve the problem in this paper.

All of the notations  in this paper are the same as in \cite
{Zh98b}.

If  $ k $  is not algebraically   closed, let  $ P $  be algebraic
closure of  $ k $. we can construct a Lie algebra  $ L_{P}=P \otimes
L $  over  $ P $, as in \cite [Section 8]{Ja62}.

By \cite [P11--14]{Ja62}, we have that

\begin {Lemma} \label {1.1}
Let  $ L $  be a vector space over  $ k $. Then  $ L $  is a Lie
algebra over field $k$ with dim  $ L $ =3 and dim  $ L^{'}=2  $  iff
there is a basis  $ e_{1} $,  $ e_{2} $,  $ e_{3} $  in  $ L $ such
that  $ [e_{1}, e_{2}] $ =0,  $ [e_{1}, e_{3}]=\alpha e_{1}+\beta
e_{2} $,  $ [e_{2}, e_{3}]=\gamma e_{1}+\delta  e_{2} $, where  $
\alpha, \beta, \gamma, \delta  \in k $, and  $ \alpha \delta -\beta
\gamma \neq 0 $.
\end {Lemma}

In this paper, we only study the Lie algebra  $ L $  in Lemma \ref
{1.1}. Set  $ A $ = $
      \left(\begin{array}{cc}
      \alpha  & \gamma  \\
      \beta  & \delta
      \end{array}\right ).
      $
Thus  $ A $  is similar to  $
                          \left(\begin{array}{cc}
                          \lambda_{1}& 0\\
                          0 & \lambda_{2}
                          \end{array}\right )
                         $
or $
     \left(\begin{array}{cc}
     \lambda_{1}& 0\\
     1& \lambda_{1}
      \end{array}\right )
      $
in the algebraic closure  $ P $  of  $ k $. Therefore, there is an
invertible matrix  $ D $  over  $ P $  such that  $ AD=D
                                                \left(\begin{array}{cc}
                                                 \lambda_{1}& 0\\
                                                  0 & \lambda_{2}
                                                 \end{array}\right )
                                                 $, or
 $ AD=D
                                                \left(\begin{array}{cc}
                                                 \lambda_{1}& 0\\
                                                  1& \lambda_{2}
                                                 \end{array}\right )
                                                 $.
                                                Let  $ Q=
                                                                 \left(\begin{array}{cc}
                                                                 D& 0\\
                                                                 0 &
                                                                \frac {1} { \lambda_{1}}
                                                                 \end{array}\right )
                                                               $
and  $ (e_{1}^{'}, e_{2}^{'}, e_{3}^{'})=(e_{1}, e_{2}, e_{3})Q $.
By computation, we have that  $ [e_{1}^{'}, e_{2}^{'}]=0 $,
 $ [e_{1}^{'}, e_{3}^{'}]=e_{1}^{'}+\beta ^{'}e_{2}^{'} $,  $ [e_{2}^{'},
e_{3}^{'}]=\delta ^{'}e_{2}^{'} $,  where  $ \beta ^{'}=0 $  and
 $ \delta ^{'}=\frac { \lambda_{1}} {\lambda_{2}} $  when  $ A $  is similar
to  $
                                       \left(\begin{array}{cc}
                                       \lambda_{1}& 0\\
                                        0 & \lambda_{2}
                                        \end{array}\right )
                                       $;
 $ \beta ^{'}=\frac {1} {\lambda_{1}} $  and  $ \delta ^{'}=1 $  when  $ A $  is
similar to  $
              \left(\begin{array}{cc}
              \lambda_{1}& 0\\
               1& \lambda_{1}
               \end{array}\right )
               $.
Let  $ Q=(q_{ij}{})_{3\times 3} $  and  $ Q^{-1}=(\bar
q_{ij})_{3\times 3} $. If  $ r =\sum_{i, j=1} ^{3}
k_{ij}(e_{i}\otimes e_{j})=\sum_{i, j=1} ^{3}
k_{ij}^{'}(e_{i}^{'}\otimes  e_{j}^{'}) $, where  $ k_{ij}\in k $, $
k_{ij}^{'}\in  P $  for  $ i, j $ =1, 2, 3, then  $$
k_{ij}^{'}=\sum_{m, n}^{3}k_{mn}\overline{q}_{im}\overline{q}_{jn}
 \hbox { and }
  k_{ij}=\sum_{m, n}^{3}k_{mn}^{'}q_{im}q_{jn} $$  for  $ i, j=1, 2, 3 $.
Obviously,  $ k_{33}=k_{33}^{'} $.

\begin {Lemma} \label{1.2}
(i)  $ k_{i3}=k_{3i} $, for  $ i $ =1, 2, 3 iff
 $ k_{i3}^{'}=k_{3i}^{'} $  for $ i $ = 1, 2, 3;

(ii)  $ k_{i3}=-k_{3i} $  for  $ i $ =1, 2, 3 iff
 $ k_{i3}^{'}=- k_{3i}^{'} $  for  $ i $ =1, 2, 3;

(iii)  $ k_{ij}=k_{ji} $  for  $ i, j $ =1, 2, 3 iff
 $ k_{ij}^{'}=k_{ji}^{'} $  for  $ i, j=1, 2, 3 $;

(iv)  $ k_{ij}=-k_{ji} $  for  $ i, j $ =1, 2, 3 iff
 $ k_{ij}^{'}=-k_{ji}^{'} $  for  $ i, j $ =1, 2, 3.

\end {Lemma}
{\bf Proof } (i) If  $ k_{i3}^{'}=k_{3i}^{'} $  for  $ i $ =1, 2, 3,
we see that
\begin {eqnarray*}
  k_{i3}  &=& \sum_{m, n}^{3} k_{mn}^{'}q_{im}q_{3n} \\
        &=& \sum_{m}^{3}k_{m3}^{'}q_{im}q_{33} \ \ ( \hbox {    since }  q_{31}=q_{32}=0
        )\\
       & =& \sum_{m}^{3}k_{3m}^{'}q_{im}q_{33}  \ \ ( \hbox {    by
       assumption})\\
       & = &  \sum_{m}^{3}k_{3m}^{'}q_{33}q_{im} \\
        &=& \sum_{m, n}^{3}k_{mn}^{'}q_{3n}q_{im} \\
        &=&  k_{3i}.
        \end {eqnarray*}
Therefore,  $ k_{i3}=k_{3i} $  for  $ i $ =1, 2, 3. The others can
be proved similarly. $\Box$

\section {The solutions of CYBE with char  $ k\neq 2  $ }

In this section, we find the general solution of CYBE for Lie
algebra  $ L $  with dim  $ L $ =3 and dim  $ L'= 2 $, where char $
k\neq 2 $.

\begin{Theorem} \label {2.1}
Let  $ L $  be a Lie algebra with a basis  $ {e_{1}, e_{2}, e_{3}} $
such that $ [e_{1}, e_{2}] $ =0,  $ [e_{1}, e_{3}]=\alpha
e_{1}+\beta e_{2} $,  $ [e_{2}, e_{3}]=\gamma  e_{1}+\delta e_{2} $,
where  $ \alpha, \beta, \gamma, \delta  \in  k $, and  $ \alpha
\delta -\beta \gamma \neq 0 $. Let $ p, q, s, t, u, v, x, y, z \in k
$. Then  $ r $  is a solution of CYBE iff  $ r $  is strongly
symmetric, or  $ r=p(e_{1}\otimes e_{2})+q(e_{2}\otimes
e_{1})+s(e_{1}\otimes e_{3})-s(e_{3}\otimes e_{1})+u(e_{2}\otimes
e_{3})-u(e_{3}\otimes e_{2})+x(e_{1}\otimes e_{1})+y(e_{2}\otimes
e_{2}) $  with  $ s(2\alpha x+\gamma (p+q))=u(2\delta y+\beta
(q+p))=u(2\alpha x+\gamma (q+p))=s(2\delta y+\beta (q+p))=(\alpha
-\delta )us+\gamma u^{2}-\beta s^{2} =s(2\gamma y+2\beta x+(\alpha
+\delta )(q+p))=u(2\gamma y+2\beta x+(\alpha +\delta )(p+q))=0 $.
\end{Theorem}

{\bf  Proof } Let  $ r =\sum_{i, j=1}^{3}k_{ij}(e_{i}\otimes
e_{j})\in  L \otimes L $, and  $ k_{ij} \in  k$, with  $ i, j $ =1,
2, 3. By computation, for all  $ i, j, n $ =1, 2, 3, we have that
the cofficient of  $ e_{j}\otimes  e_{i}\otimes  e_{i} $  in
 $ [r ^{12}, r ^{23}] $  is zero and the cofficient of  $ e_{i}\otimes  e_{i}\otimes  e_{j} $  in
 $ [r ^{13}, r ^{23}] $  is zero.

We can obtain the following equations by seeing the cofficient of
 $ e_{i} \otimes  e_{j} \otimes  e_{n} $  in  $ [r ^{12}, r
^{13}]+[r ^{12}, r ^{23}]+[r ^{13}, r ^{23}] $, as in \cite
[Proposition 2.6]{Zh98b}. To simplify notation, let
 $ k_{11}=x, k_{22}=y, k_{33}=z, k_{12}=p, k_{21}=q, k_{13}=s,
k_{31}=t, k_{23}=u, k_{32}=v $.

  (1)  $ -\alpha sx+\alpha  xt-\gamma sq+\gamma pt=0 $;

  (2)  $ -\beta up+\beta qv-\delta uy+\delta yv=0 $;

  (3) $ -\alpha vs+\alpha pz-\gamma uv+\gamma yz-\beta s^{2}+\beta xz-\delta su+\delta zp=0 $;

  (4) $ -\beta xz+\beta st-\delta zq+\delta ut-\alpha ut+\alpha qz-\gamma uv+\gamma yz=0 $;

  (5) $ -\alpha zp+\alpha tv-\gamma zy+\gamma uv-\beta zx+\beta ts-\delta zp+\delta vs=0 $;

  (6) $ -\alpha zp+\alpha sv-\gamma zy+\gamma uv-\beta st+\beta xz-\delta sv+\delta zp=0
  $;

  (7) $ -\beta zx+\beta tt-\delta zq+\delta vt-\alpha zq+\alpha tu-\gamma zy+\gamma uv=0 $;

  (8) $ -\beta st+\beta xz-\delta tu+\delta qz-\alpha us+\alpha qz-\gamma uu+\gamma yz=0 $;

  (9) $ -\alpha tp+\alpha xv-\gamma ty+\gamma qv-\alpha sp+\alpha vx-\gamma sy+\gamma pv=0 $;

  (10) $ -\alpha ux+\alpha qt-\gamma uq+\gamma yt-\alpha ux+\alpha qs-\gamma up+\gamma ys=0 $;

  (11) $ -\alpha st+\alpha xz-\gamma tu+\gamma qz-\alpha s^{2}+\alpha xz-\gamma su+ypz=0 $;

  (12) $ -\alpha zx+\alpha tt-\gamma zq+\gamma vt-\alpha zx+\alpha st-\gamma zp+\gamma vs=0 $;

  (13) $ -\beta vx+\beta pt-\delta vq+\delta yt-\beta ux+\beta qt-\delta uq+\delta yt=0 $;

  (14) $ -\beta sp-\beta xv-\delta sy+\delta pv-\beta sq+\beta xu-\delta sy+\delta pu=0 $;

  (15) $ -\beta vs+\beta pz-\delta vu+\delta yz-\beta su+\beta zq-\delta u^{2}+\delta yz=0 $;

  (16) $ -\beta zp+\beta tv-\delta zy+\delta vv-\beta zq+\beta tu-\delta zy+\delta vu=0 $;

  (17) $ -\gamma zq+\gamma ut-\gamma sv+\gamma pz=0 $;

  (18) $ -\alpha vx+\alpha pt-\gamma vq+\gamma yt-\beta sx+\beta xt-\delta sq+\delta pt-\alpha sq+\alpha xu-\gamma sy+\gamma pu=0 $;

  (19) $ -\beta pt+\beta vx-\delta ty+\delta qv-\alpha up+\alpha qv-\gamma uy+\gamma yv-\beta ux+\beta qs-\delta up+\delta ys=0 $;

  (20) $ -\beta zp+\beta sv-\delta zy+\delta uv-\beta ut+\beta qz-\delta uv+\delta yz=0 $;

  (21) $ -\beta zs+\beta tz-\delta zu+\delta vz=0 $;

  (22) $ -\alpha zs+\alpha tz-\gamma zu+\gamma vz=0 $;

It is clear that  $ r $  is a solution of CYBE iff (1)-(22) hold.

By simple computation, we have the sufficiency. Now we show the
necessity, If  $ k_{33}\neq 0 $, then  $ k_{33}^{'}\neq 0' $  and so
$ r $  is a strongly symmetric element in  $ L_{P}\otimes L_{P} $ by
\cite [The proof of Proposition 1.6 ] {Zh98b}. Thus  $ r $  is a
strongly symmetric element in  $ L\otimes L $. If  $ k_{33} $ =0,
then  $ k_{33}^{'} $ =0. By \cite [Proposition 1.6] {Zh98b}, we have
that  $ k_{i3}^{'}=-k_{3i}^{'} $  for  $ i $ =1, 2, 3, which implies
that  $ k_{i3}=-k_{3i} $  for  $ i $ =1, 2, 3 by Lemma \ref {1.2}.

It immediately follows from (1)-(22) that

  (23) $ s(-2\alpha x-\gamma (q+p))=0 $;

  (24) $ u(2\delta y+\beta (q+p))=0 $;

  (25) $ \gamma u^{2}-\beta s^{2}+(\alpha -\delta )us=0 $;

  (26) $ u(2\alpha x+\gamma (q+p))=0 $;

  (27) $ s(2\delta y+\beta (q+p))=0 $;

  (28) $ 2\alpha ux-2\gamma ys-2\beta xs+(-s(\alpha +\delta )+\gamma u)(q+p)=0 $;

  (29) $ -2\beta ux+2\delta ys-2\gamma uy+(-u(\alpha +\delta )+\beta s)(q+p)=0
  $;\\
By (26), (27), (28) and (29), we have that

  (30)  $ s(2\gamma y+2\beta x+(\alpha +\delta )(q+p))=0 $;

  (31)  $ u(2\gamma y+2\beta x+(\alpha +\delta )(p+q))=0 $. $\Box$

\begin{Example} \label {2.2}
Let  $ L $  be a Lie algebra over real field  $ R $  with dim  $ L $
=3 and dim  $ L^{'}=2. $  If there is complex characteristic root
 $ \lambda_{1}=a+bi $  of  $ A $  and the root is not real, then  $ r\in
L\otimes L $  is a solution of CYBE iff  $ r $  is strongly
symmetric, or  $ r=p(e_{1}\otimes e_{2})+q(e_{2}\otimes
e_{1})+x(e_{1}\otimes e_{1})+y(e_{2}\otimes e_{2}) $  for any  $ p,
q, x, y \in R $.
\end{Example}

{\bf Proof.} There are two different characteristic roots:
 $ \lambda_{1}=a+bi $  and  $ \lambda_{2}=a-bi $, where  $ a, b\in  R $, Thus
 $ A $  must be similar to $
       \left(\begin{array}{cc}
       a& b\\
         -b& a
         \end{array}\right )
        $.  By Theorem \ref {2.1}, we can complete the proof. $\Box$

\section{The solutions of CYBE with char  $ k $ =2}

In this section, we find the general solution of CYBE for Lie
algebra  $ L $  with dim  $ L $ =3 and dim  $ L^{'}=2 $, where char
 $ k $ =2.

\begin{Theorem} \label{3.1}
Let  $ L $  be a Lie algebra with a basis  $ {e_{1}, e_{2}, e_{3}} $
such that  $ [e_{1}, e_{2}] $ =0,  $ [e_{1}, e_{3}]=\alpha
e_{1}+\beta e_{2} $,  $ [e_{2}, e_{3}]=\gamma e_{1}+\delta e_{2} $,
where  $ \alpha, \beta, \gamma, \delta  \in  k $, and  $ \alpha
\delta -\beta \gamma \neq 0 $. Let $ p, q, s, t, u, v, x, y, z \in k
$  and  $ A=
       \left(\begin{array}{cc}
       \alpha & \gamma \\
         \beta & \delta
         \end{array}\right )
      $.

 (I) If two characteristic roots of  $ A $  are equal and  $ A $  is
 similar to a diagonal matrix in the algebraic closure  $ P $  of
  $ k $, then  $ r $  is a solution of CYBE in  $ L $  for any  $ r \in L\otimes L $;

 (II) If the condition in Part (1)does not hold, then  $ r $  is a
 solution of CYBE in  $ L $  iff  $ r=p(e_{1}\otimes  e_{1})+p(e_{2}\otimes
e_{1})+s(e_{1}\otimes  e_{3})+s(e_{3}\otimes  e_{1})+u(e_{2}\otimes
e_{3})-+u(e_{3}\otimes  e_{2})+x(e_{1}\otimes  e_{1})+y(e_{2}\otimes
e_{2})+z(e_{3}\otimes  e_{3}) $  with $ \alpha us+\alpha pz+\gamma
u^{2}+\gamma yz+\beta s^{2}+\beta xz+\delta su+\delta zp=0 $  and
 $ z\neq 0 $; or $ r=p(e_{1}\otimes  e_{2})+q(e_{2}\otimes
e_{1})+s(e_{1}\otimes  e_{3})+s(e_{3}\otimes  e_{1})+u(e_{2}\otimes
e_{3})+u(e_{3}\otimes  e_{2})+x(e_{1}\otimes  e_{1})+y(e_{2}\otimes
e_{2})$  with  $ s\gamma (p+q)=u\beta (p+q)=u\gamma (p+q)=s\beta
(p+q)=(\alpha +\delta )us+\gamma u^{2}+\beta s^{2} $  = $ s(\alpha
+\delta )(p+q)=u(\alpha +\delta )(p+q)=0 $.

\end{Theorem}

{\bf Proof}
 We only show the necessity since the sufficiency  can easily be shown.
By the proof of Theorem  \ref {2.1}, there exists an invertible
matrix  $ Q $  such that  $ (e_{1}^{'}, e_{2}^{'},
e_{3}^{'})=(e_{1}, e_{2}, e_{3})Q $  and  $ [e_{1}^{'}, e_{2}^{'}]=0
$,  $ [e_{1}^{'}, e_{3}^{'}]=e_{1}^{'}+\beta ^{'}e_{2}^{'} $,
 $ [e_{2}^{'}, e_{3}^{'}]=\delta ^{'}e_{2}^{'} $. We use the notations
before Lemma \ref {1.2}.  By [11, Proposition2.4],
 $ k_{i3}^{'}=k_{3i}^{'} $  for  $ i $ =1, 2, 3, which implies that
 $ k_{3i}=k_{i3} $  for  $ i $ =1, 2, 3 by Lemma \ref {1.2}.

(I) If  $ A $  is similar to   $
                            \left(\begin{array}{cc}
                              \lambda_{1}& 0\\
                                0& \lambda_{1}
                              \end{array}\right )
                            $, then  $ \beta ^{'}=0 $  and
  $ \delta'=1 $. By [11, Proposition 2.4], we have Part (I).

 (II) Let  $ A $  be not similar to   $
                                  \left(\begin{array}{cc}
                                    \lambda_{1}& 0\\
                                     0& \lambda_{1}
                                   \end{array}\right )
                                    $.

  (a). If  $ z\neq 0 $,  then  $ k_{33}^{'}\neq 0 $ \. Thus
   $ k_{12}^{'}=k_{21}^{'} $  by \cite [Proposition 2.4]{Zh99f}, which implies
   $ k_{12}=k_{21} $. It is straightforward to check that relation
  (1)-(22) in the proof of Theorem \ref {2.1} hold iff
   $ \alpha us+\alpha pz+\gamma u^{2}+\gamma yz+\beta s^{2}+\beta xz+\delta su+\delta zp=0. $

  (b). If  $ z $ =0, then we can obtain that  $ r $  is the second case
  in Part (II) by using the method similar to the proof of
  Theorem \ref {2.1}. $\Box$

  \section{Coboundary Lie bialgebras}

In this section, using the general solution, which are obtained in
the section above, of CYBE in Lie algebra  $ L $  with dim $ L $ =3
and dim $ L^{'} $ =2, we give the sufficient and necessary
conditions which ( $ L, [\ ], \Delta_{r}, r $ ) is a coboundary (or
triangular)Lie bialgebra.

\begin{Theorem} \label {4.1}
Let  $ L $  be a Lie algebra with a basis  $ {e_{1}, e_{2}, e_{3}} $
such that $ [e_{1}, e_{2}] $ =0,  $ [e_{1}, e_{3}]=\alpha
e_{1}+\beta e_{2} $,  $ [e_{2}, e_{3}]=\gamma e_{1}+\delta e_{2} $,
where  $ \alpha, \beta, \gamma, \delta  \in  k $, and  $ \alpha
\delta -\beta \gamma \neq 0 $. Set
  $ A=
       \left(\begin{array}{cc}
      \alpha & \gamma \\
         \beta & \delta
         \end{array}\right )
      $. Let  $ p, u, s\in  k $,
      $ r\in  Im(1-\tau) $  and  $ r=p(e_{1}\otimes  e_{2})-p(e_{2}\otimes
     e_{1})+s(e_{1}\otimes  e_{3})-s(e_{3}\otimes  e_{1})+u(e_{2}\otimes  e_{3})-u(e_{3}\otimes
     e_{2}) $.
Then

(I) \ ( $ L,  [\ ], \Delta_{r},  r $ ) is a coboundary  Lie
bialgebra iff  $$ (s, u)
                       \left(\begin{array}{cc}
                       \beta \delta +\alpha \beta & -\beta \gamma -\alpha ^{2}\\
                         \delta ^{2}+\gamma \beta & -\delta \gamma -\gamma \alpha
                        \end{array}\right )
                                                   \left(\begin{array}{cc}
                             s\\
                             u
                             \end{array}\right )
                            =0; $$

(II) If two characteristic roots of  $ A $  are equal and  $ A $  is
similar to a diagonal matrix in the algebraic closure  $ P $  of  $
k $  with char $ k $ =2,  then ( $ L, [\ ], \Delta _ r,  r $ ) is a
triangular Lie bialgebra for any  $ r\in  Im(1-\tau ) $;

(III) If the condition in Part(II) does not hold, then( $ L, [\ ],
\Delta _r,  r $ ) is a triangular Lie bialgebra iff  $ -\beta
s^{2}+\gamma u^{2}+(\alpha -\delta )us=0 $.
\end {Theorem}

{\bf Proof } We can complete the proof as in the proof of \cite
[Theorem 3.3]{Zh98b}. $\Box$

\begin{Example} \label{4.2}

Under Example \ref {2.2}, and  $ r\in  Im(1-\tau) $, we have the
following :

(I) ( $ L, [\ ], \Delta_r,  r $ ) is a coboundary Lie bialgebra iff
 $ r=p(e_{1}\otimes  e_{2})-p(e_{2}\otimes  e_{1})+s(e_{1}\otimes
e_{3})-s(e_{3}\otimes  e_{1})+u(e_{2}\otimes  e_{3})-u(e_{3}\otimes
e_{2}) $  with  $ a(s^{2}+u^{2})=0 $;

(II) ( $ L, [\ ], \Delta_r, r $ )is a triangular Lie bialgebra iff
 $ r=p(e_{1}\otimes  e_{2})-p(e_{2}\otimes  e_{1}) $.

\end {Example}

\end {document}